\documentclass[twocolumn]{autart}

\usepackage{epsfig} 
\usepackage{fancyhdr}

\usepackage{amsfonts}
\usepackage{mathrsfs}
\usepackage{colortbl}
\usepackage{amssymb,amsmath,graphicx}
\usepackage{caption}

\usepackage{subfigure}
\usepackage{float}

\newtheorem{Example}{Example}
\newtheorem{Remark}{Remark}[section]

\newtheorem{Problem}{Problem}
\newenvironment{Proof}{\noindent{\em Proof:\/}}{\hfill $\Box$\par}
\newtheorem{Theorem}{Theorem}[section]
\newtheorem{Lemma}{Lemma}[section]

\newtheorem{Assumption}{Assumption}

\newcommand{\EQQ}{\begin{eqnarray*}}
\newcommand{\ENN}{\end{eqnarray*}}
\newcommand{\EQ}{\begin{eqnarray}}
\newcommand{\EN}{\end{eqnarray}}

\begin{document}

\begin{frontmatter}

\title{Cooperative Global Robust Stabilization for a Class of Nonlinear Multi-Agent Systems and its Application
\thanksref{footnoteinfo}} 

\thanks[footnoteinfo]{ This work has been supported in part by the Research Grants
Council of the Hong Kong Special Administration Region under grant
No. 412813, and in part by National Natural Science Foundation of
China under Project 61174049. Corresponding author: Jie Huang. Tel.
+852-39438473. Fax +852-39436002.}

\author[Liu-Huang]{Wei Liu}\ead{wliu@mae.cuhk.edu.hk},
\author[Liu-Huang]{Jie Huang}\ead{jhuang@mae.cuhk.edu.hk}.

\address[Liu-Huang]{Shenzhen Research Institute,  The Chinese University of Hong Kong, and Department of Mechanical and Automation
Engineering, The Chinese University of Hong Kong.}

\begin{keyword}
Cooperative control, global robust stabilization, multi-agent systems, nonlinear systems,  switched control.
\end{keyword}

\begin{abstract}
This paper studies the cooperative global robust stabilization problem for a class of nonlinear multi-agent systems.
The problem is motivated from the study of the cooperative global robust output regulation problem for the class of nonlinear multi-agent systems in normal form with unity relative degree
which was studied recently under the conditions that the switching network is undirected and some nonlinear functions satisfy certain growth condition.
 We first solve the stabilization problem by using the multiple Lyapunov functions approach and the average dwell time method. Then, we apply this result to the cooperative global robust output regulation problem for the class of nonlinear systems in normal form with unity relative degree under directed switching network, and have removed the conditions that the switching network is undirected and some nonlinear functions satisfy certain growth condition.
\end{abstract}

\end{frontmatter}

\section{Introduction}
Consider a class of cascade-connected multi-agent nonlinear systems as follows:
\begin{equation}\label{system1a}
\begin{split}
  &\dot{Z}_{i} = F_{i}(Z_{i},e_{i},d(t))\\
  &\dot{e}_{i} = G_{i}(Z_{i},e_{i},d(t))+b_{i} \bar{u}_{i}, ~i=1,\cdots,N \\
\end{split}
\end{equation}
where $(Z_{i},e_{i})\in\mathbb{R}^{n_{i}}\times\mathbb{R}$ is the state, $\bar{u}_{i}\in \mathbb{R}$ is the input, $e_{i}$ is the output, $b_{i}$ is an unknown real number, and $d:[t_{0},\infty){\rightarrow}\mathbb{D}\subset\mathbb{R}^{n_{d}}$ with $\mathbb{D}$ a known compact subset represents  external disturbance and/or  parameter variations.
It is assumed that $b_{m} \leq |b_i| \leq b_{M}$ for some known positive real numbers $b_{m}$ and $b_{M}$, and
the functions  $F_{i}:\mathbb{R}^{n_{i}}\times\mathbb{R}\times\mathbb{R}^{n_{d}}{\rightarrow}\mathbb{R}^{n_{i}}$, and $G_{i}:\mathbb{R}^{n_{i}}\times\mathbb{R}\times\mathbb{R}^{n_{d}}{\rightarrow}\mathbb{R}$ are both sufficiently smooth and
 satisfy $F_{i}(0,0,d(t))=0$ and $G_{i}(0,0,d(t))=0$ for all $d(t)\in\mathbb{R}^{n_{d}}$.

To describe our control law, let $\sigma:[0,\infty)\rightarrow\mathcal {P} = \{1,2,
\cdots, n_{0} \}$ for some integer $n_{0}  > 0$, be  a
piecewise constant switching signal, and
$S_{ave}[\tau_{d},N_{0}]$ be the set of all signals possessing
the property of average dwell-time
$\tau_{d}$ with chatter bound $N_{0}$ \cite{Liberzon1,Liu2}, $H_{p} =[h_{ij}^{p}]\in \mathbb{R}^{N \times N}$, $ p = 1, \cdots, n_0$, be some $\mathcal{M}$ matrices\footnote{A matrix $M\in\mathbb{R}^{N\times N}$ is called an $\mathcal{M}$ matrix if all of its non-diagonal elements are non-positive and all of its eigenvalues have positive real parts.}.
Define a piecewise switching matrix $H_{\sigma (t)}  \in  {R}^{N \times N}$ such that, over each interval $[t_i, t_{i+1})$, $ H_{\sigma (t)}  = H_{p} $ for
some integer $1 \leq p \leq n_{0}$. Denote the elements of  $H_{\sigma (t)}$ by $ h_{ij} (t)$,  $i,j = 1, \cdots, N$, and define the virtual output of (\ref{system1a}) as
\begin{equation}\label{evi}
\begin{split}
e_{v\sigma(t)i} =  \sum_{j=1}^N h_{ij} (t) e_j,\ i=1,\cdots,N.
\end{split}
\end{equation}
Then,  we describe our control law as follows:
\begin{equation}\label{ui2a}
\begin{split}
\bar{u}_{i}=k_{i}(e_{v\sigma(t)i}), \;\;
i=1,...,N
\end{split}
\end{equation}
where the functions $k_{i}(\cdot)$, $i=1,\cdots,N$, are sufficiently smooth vanishing at the origin.
 Such a control law is called a
distributed switched output feedback control law, since we can only use $e_{v\sigma(t)i}$ instead of $e_{i}$ for feedback control due to the communication constraints, which will be further elaborated in Section \ref{Application}.

%
%
%
%

 \begin{Problem}\label{Definition1a}
Given the multi-agent system \eqref{system1a}, a set of $\mathcal{M}$ matrices $H_p \in \mathbb{R}^{N \times N}$, $ p = 1, \cdots, n_0$,  and some compact subset
$\mathbb{D}\subset\mathbb{R}^{n_{d}}$ with $0\in\mathbb{D}$, find $\tau_{d},N_{0}$, and  a control law of the
form (\ref{ui2a})  such that, for any $d(t)\in \mathbb{D}$,  and any $ \sigma \in S_{ave}[\tau_{d},N_{0}]$, the equilibrium point of the closed loop system composed of \eqref{system1a} and \eqref{ui2a} at the origin is
globally asymptotically stable. \end{Problem}


The above cooperative global robust stabilization problem is of interest on its own and has not been studied before, since \eqref{ui2a} is a switched control law which results in a switched closed loop system. On the other hand, it is  motivated from the study of the cooperative global robust output regulation problem for a class of nonlinear multi-agent systems with switching network.
In fact, we will show that our main result is applicable to the cooperative global robust output regulation problem for the unity relative degree nonlinear multi-agent systems with switching network in Section \ref{Application}. It is also noted that the result in \ref{Application} includes some existing results as special cases.
If $n_0 = 1$, then  $H_{\sigma(t)}$ reduces to a constant matrix. For this special case, the above problem has been studied in \cite{Dongyi1}
under the assumption that $H_1$ is a symmetric $\mathcal{M}$ matrix. More recently, for the case where $n_0$ is any positive integer,
the above problem was studied in \cite{Liu2}  under the assumption that $H_p$, $p = 1, \cdots, n_{0}$,  are all symmetric $\mathcal{M}$ matrices and the nonlinear function $G_i$ satisfy certain
growth condition  \footnote{See Assumption \ref{Ass3.3o} for the definition of growth condition.}. In this paper, we will further study the above problem for the more general case where
 the nonlinear function $G_i$ does not satisfy any growth condition, and $H_p$, $p = 1, \cdots, n_{0}$, are any  $\mathcal{M}$ matrices, i.e., do not need to be symmetric.  As an application of our main result, we will obtain the solution of the same problem studied in \cite{Liu2}
without assuming that the nonlinear functions $G_i$  satisfy any
growth condition, and the communication graph of the multi-agent system is undirected for all $ t \geq 0$.

Without assuming the growth condition, and the undirectedness of the communication graph,   the problem in this paper is  technically  much more challenging than the one in \cite{Liu2}.
To overcome these difficulties, we need to develop a changing supply pair technique for exponentially input-to-state stable (exp-ISS)  nonlinear systems  to remove the growth condition of the nonlinear function $G_i$,
and we need to employ a non-quadratic function for the closed-loop system in Section \ref{Result} to handle  the directed communication graph.


{\bf Notation.} For any column vectors $a_i$, $i=1,...,s$, denote $\mbox{col}(a_1,...,a_s)=[a_1^T,...,a_s^T]^T$.
A function $\bar{\rho}(\cdot)$ is called $\mathcal{SN}$ function if it is a smooth non-decreasing function $\bar{\rho}:[0,\infty){\rightarrow}[0,\infty)$ satisfying $\bar{\rho}(s)>0$ for all $s>0$.
  $\bar{\alpha}(s)\!=\!O(\alpha(s))$ as $s\rightarrow0^{+}$ means $\lim\limits_{s\rightarrow0^{+}} \sup\frac{\bar{\alpha}(s)}{\alpha(s)}\!<\!\infty$. The notation $\lambda_{\min}(A)$ denotes the minimum eigenvalue of a symmetric real matrix $A$.

%

\section{Preliminaries}\label{Preliminary}
In this section, we will establish a technical lemma. Consider a general nonlinear system:
 \begin{equation}\label{system0a}
\begin{split}
\dot{x}=f(x,u)
\end{split}
\end{equation}
where $x\in\mathbb{R}^{n}$ is the state, $u\in\mathbb{R}^{m}$ is the input, $f:\mathbb{R}^{n}\times\mathbb{R}^{m}{\rightarrow}\mathbb{R}^{n}$ is locally Lipschitz and $f(0,0)=0$.
It is known from \cite{Sontag1} that system (\ref{system0a}) is said to be input-to-state stable (ISS) if there exists a $C^1$ function $V:  {R}^{n}{\rightarrow} [0,\infty)$  such that, for all $x,u$,
  \begin{equation}\label{dotVx1b}
\begin{split}
&\alpha_{1}(\|x\|)  \leq V(x)\leq \alpha_{2}(\|x\|) \\
&\frac{\partial V}{\partial x}f (x,u) \leq-\alpha(\|x\|) +\gamma(\|u\|)
\end{split}
\end{equation}
for some class $\mathcal{K}_{\infty}$ functions
 $\alpha_{1}(\cdot)$, $\alpha_{2}(\cdot)$, $\alpha(\cdot)$ and some class $\mathcal{K}$ function $\gamma(\cdot)$. The pair of  functions $(\gamma,\alpha)$ is called a supply pair for system (\ref{system0a}), and $V (x)$ is called an ISS Lyapunov function of (\ref{system0a}).
Moreover,  if system (\ref{system0a}) is ISS,  then, for any class $\mathcal{K}_{\infty}$ function $\bar{\alpha}(\cdot)$ satisfying $\bar{\alpha}(s)\!=\!O(\alpha(s))$ as $s\rightarrow0^{+}$,
 there exists $C^1$ function $\bar{V}(x)$ such that,  for all $x,u$,
   \begin{equation}\label{dotVx1a}
\begin{split}
& \bar{\alpha}_{1}(\|x\|)\leq \bar{V}(x)\leq \bar{\alpha}_{2}(\|x\|) \\
&\frac{\partial \bar{V}}{\partial x}f (x,u) \leq- \bar{\alpha}(\|x\|) +\bar{\gamma}(\|u\|)
\end{split}
\end{equation}
for some class $\mathcal{K}_{\infty}$ functions
 $\bar{\alpha}_{1}(\cdot)$, $\bar{\alpha}_{2}(\cdot)$ and some class $\mathcal{K}$ function $\bar{\gamma}(\cdot)$ \cite{Sontag1}.  This result is called the changing supply pair technique which plays a key role in finding a suitable Lyapunov function
for a nonlinear system to conclude the asymptotic stability of its origin. However, as will be explained in Section III,  this version of the changing supply pair technique is not adequate for
handling the stability of switched systems. We need to further establish a lemma for the following class of nonlinear systems:
 \begin{equation}\label{system2a}
\begin{split}
\dot{x}=f(x,u,d(t))
\end{split}
\end{equation}
where $x\in\mathbb{R}^{n}$ is the state, $u\in\mathbb{R}^{m}$ is the input, $d:[t_{0},\infty){\rightarrow}\mathbb{D}\subset\mathbb{R}^{n_{d}}$ with $\mathbb{D}$ some non-empty set, represents external unpredictable disturbance and/or internal parameter variation, $f:\mathbb{R}^{n}\times\mathbb{R}^{m}\times\mathbb{D}{\rightarrow}\mathbb{R}^{n}$ is piecewise continuous in $d(t)$ and locally Lipschitz, $d(t)$ is piecewise continuous in $t$ and $f(0,0,d(t))=0$ for any $d(t)\in\mathbb{R}^{n_{d}}$.

 \begin{Lemma}\label{Lemma0}
Suppose that there exists a $C^{1}$ function $V(x,t)$ 
such that, for any $d(t)\in\mathbb{D}$,
 \begin{equation}\label{Vx1}
\begin{split}
\alpha_{1}(\|x\|)\leq V(x,t)\leq \alpha_{2}(\|x\|),\ \ \forall\ x
\end{split}
\end{equation}
 \begin{equation}\label{dotVx1}
\begin{split}
\dot{V}(x,t)\leq
-\lambda V(x,t)+\beta(u),\ \ \forall\ x,u
\end{split}
\end{equation}
where $\alpha_{1}(\cdot)$ and $\alpha_{2}(\cdot)$ are some known class $\mathcal{K}_{\infty}$ functions,  $\lambda$ is some known positive real number and $\beta(\cdot)$ is some known smooth positive definite function. Then,

(i) For any function $\bar{\rho}\in\mathcal{SN}$, the following $C^1$ function \begin{equation}\label{bVx2}
\begin{split}
\bar{V}(x,t)=\int_{0}^{V(x,t)}\bar{\rho}(s)ds,
\end{split}
\end{equation}
is an ISS Lyapunov function in the sense that
 \begin{equation}\label{bVx1}
\begin{split}
\bar{\alpha}_{1}(\|x\|)\leq \bar{V}(x,t)\leq \bar{\alpha}_{2}(\|x\|),\ \ \forall\ x
\end{split}
\end{equation}
 \begin{equation}\label{dotbVx1}
\begin{split}
\dot{\bar{V}}(x,t)\leq
- \dfrac{\lambda}{2}  \bar{V}(x,t)+\bar{\beta}(u),\ \ \forall\ x,u, \;  \forall d(t)\in\mathbb{D}
\end{split}
\end{equation}
 for some class $\mathcal{K}_{\infty}$ functions $\bar{\alpha}_{1}(\cdot)$ and $\bar{\alpha}_{2}(\cdot)$,  some smooth positive definite function $\bar{\beta}(\cdot)$.

(ii) For any class $\mathcal{K}_{\infty}$ function $\bar{\alpha}(\cdot)$ satisfying  $\bar{\alpha}(s)=O(\alpha_{1}(s))$ as $s\rightarrow0^{+}$, there exists some function $\bar{\rho}\in\mathcal{SN}$
such that, for any $d(t)\in\mathbb{D}$,  $\bar{V}(x,t)$ satisfies (\ref{bVx1}), and the following
 \begin{equation}\label{dotbVx3a}
\begin{split}
\dot{\bar{V}}(x,t) \leq&-(\bar{\lambda}\bar{V}(x,t)+\bar{c}\bar{\alpha}(\|x\|))+\bar{\beta}(u),\ \ \forall\ x,u
\end{split}
\end{equation}
for   any  positive real number $ 0< \bar{\lambda} < \dfrac{\lambda}{2}$ and positive real number $\bar{c} =  \dfrac{\lambda}{2} - \bar{\lambda}$.
\end{Lemma}

 \begin{Proof}
Part (i):  It is easy to see that
  \begin{equation}\label{bVx3}
\begin{split}
\bar{V}(x,t)\leq\bar{\rho}(V(x,t))V(x,t),
\end{split}
\end{equation}
and equations \eqref{Vx1} and \eqref{bVx2} imply that
 \begin{equation}\label{bVx4}
\begin{split}
&\bar{\alpha}_{1}(\|x\|)\leq\int_{0}^{\alpha_{1}(\|x\|)}\bar{\rho}(s)ds\leq \bar{V}(x,t)\\
&\leq\int_{0}^{\alpha_{2}(\|x\|)}\bar{\rho}(s)ds\leq \bar{\alpha}_{2}(\|x\|),\ \ \forall\ x
\end{split}
\end{equation}
for some class $\mathcal{K}_{\infty}$ functions $\bar{\alpha}_{1}(\cdot)$ and $\bar{\alpha}_{2}(\cdot)$.

We now show that along the trajectory of $\dot{x}=f(x,u,d(t))$
 \begin{equation}\label{dotbVx2}
\begin{split}
\dot{\bar{V}}(x,t)\leq&\bar{\rho}(V(x,t))\big(-\lambda V(x,t)+\beta(u)\big)\\
\leq&-\dfrac{\lambda}{2}\bar{\rho}(V(x,t))V(x,t)+\bar{\rho}(\dfrac{2}{\lambda}\beta(u))\beta(u).
\end{split}
\end{equation}
For this purpose, consider the following two cases.
\begin{enumerate}
  \item $\dfrac{\lambda}{2}V(x,t)\geq\beta(u)$: In this case,
\begin{equation*}
\begin{split}
&\bar{\rho}(V(x,t))\big(-\lambda V(x,t)+\beta(u)\big)\\
\leq&\bar{\rho}(V(x,t))\big(-\lambda V(x,t)+\dfrac{\lambda}{2} V(x,t)\big)\\
=&-\dfrac{\lambda}{2}\bar{\rho}(V(x,t))V(x,t)\\
\leq&-\dfrac{\lambda}{2}\bar{\rho}(V(x,t))V(x,t)+\bar{\rho}(\dfrac{2}{\lambda}\beta(u))\beta(u).
\end{split}
\end{equation*}
  \item $\dfrac{\lambda}{2}V(x,t)<\beta(u)$: In this case, we have
\begin{equation*}
\begin{split}
\bar{\rho}(V(x,t))\leq\bar{\rho}(\dfrac{2}{\lambda}\beta(u))\\
\end{split}
\end{equation*}
which implies that
\begin{equation*}
\begin{split}
&\bar{\rho}(V(x,t))\big(-\lambda V(x,t)+\beta(u)\big)\\
\leq&-\lambda\bar{\rho}(V(x,t))V(x,t)+\bar{\rho}(\dfrac{2}{\lambda}\beta(u))\beta(u)\\
\leq&-\dfrac{\lambda}{2}\bar{\rho}(V(x,t))V(x,t)+\bar{\rho}(\dfrac{2}{\lambda}\beta(u))\beta(u).
\end{split}
\end{equation*}
\end{enumerate}

Choose a smooth positive definite function  $\bar{\beta}(\cdot)$ such that
 \begin{equation}\label{bsigma1}
\begin{split}
\bar{\beta}(u)\geq\bar{\rho}(\dfrac{2}{\lambda}\beta(u))\beta(u)
\end{split}
\end{equation}

From \eqref{bVx3} , \eqref{dotbVx2}, and \eqref{bsigma1},   we have
 \begin{equation}\label{dotbVx3}
\begin{split}
\dot{\bar{V}}(x,t)\leq&-\dfrac{\lambda}{2}\bar{\rho}(V(x,t))V(x,t)+\bar{\beta}(u)\\
\leq&- \dfrac{\lambda}{2} \bar{V}(x,t)+\bar{\beta}(u).\\
\end{split}
\end{equation}

Part (ii):  Since $\bar{\alpha}(s)=O(\alpha_{1}(s))$ as $s\rightarrow0^{+}$, by Lemma 2 of \cite{Sontag1}, 
   it is always possible to find a function $\bar{\rho}\in\mathcal{SN}$ such that, for all $x$,
 \begin{equation}\label{brho0}
\begin{split}
\bar{\rho}(\alpha_{1}(\|x\|))\alpha_{1}(\|x\|)\geq\bar{\alpha}(\|x\|).
\end{split}
\end{equation}
Thus
 \begin{equation}\label{brho1}
\begin{split}
\bar{\rho}(V(x,t))V(x,t)\!\geq\!\bar{\rho}(\alpha_{1}(\|x\|))\alpha_{1}(\|x\|)\geq\bar{\alpha}(\|x\|).\\
\end{split}
\end{equation}
Choose a real number $\bar{\lambda}$ satisfying $0<\bar{\lambda}<\dfrac{\lambda}{2}$ and let $\bar{c}=\dfrac{\lambda}{2}-\bar{\lambda}$. Then, from \eqref{bVx3},
\eqref{dotbVx3} and \eqref{brho1} , we have
 \begin{equation*}
\begin{split}
\dot{\bar{V}}(x,t)\leq&-\dfrac{\lambda}{2}\bar{\rho}(V(x,t))V(x,t)+\bar{\beta}(u)\\
=&-(\bar{\lambda}+\bar{c})\bar{\rho}(V(x,t))V(x,t)+\bar{\beta}(u)\\
\leq&-(\bar{\lambda}\bar{V}(x,t)+\bar{c}\bar{\alpha}(\|x\|))+\bar{\beta}(u).\\
\end{split}
\end{equation*}
Thus the proof is completed.
 \end{Proof}

\begin{Remark}\label{Remark1} Lemma \ref{Lemma0} and its proof can be viewed as an extension of the main result in \cite{Sontag1}.
If we let $\bar{\alpha}(\cdot)$ be any smooth positive definite function $\hat{\alpha}(\cdot)$ in Lemma \ref{Lemma0}, then there exists a class $\mathcal{K}_{\infty}$ function $\tilde{\alpha}(\cdot)$ satisfying $\tilde{\alpha}(s)=O(s^2)$ as $s\rightarrow0^{+}$, and  $\tilde{\alpha}(\|x\|)\geq\hat{\alpha}(x)$ for any $x\in\mathbb{R}^{n}$. Thus,
if $\alpha_{1}(\cdot)$ satisfies $\lim\limits_{s\rightarrow0^{+}}\sup\frac{s^{2}}{\alpha_{1}(s)}<\infty$, then $\tilde{\alpha}(s)=O(\alpha_{1}(s))$ as $s\rightarrow0^{+}$.
Then we conclude that if $\alpha_{1}(\cdot)$ satisfies $\lim\limits_{s\rightarrow0^{+}}\sup\frac{s^{2}}{\alpha_{1}(s)}<\infty$, then, for any smooth positive definite function $\hat{\alpha}(\cdot)$,
by Part (ii) of Lemma \ref{Lemma0}, we have
 \begin{equation}\label{dotbVx4}
\begin{split}
\dot{\bar{V}}(x,t)\leq&-(\bar{\lambda} \bar{V}(x,t)+\bar{c}\tilde{\alpha}(\|x\|))+\bar{\beta}(u)\\
\leq&-(\bar{\lambda} \bar{V}(x,t)+\bar{c}\hat{\alpha}(x))+\bar{\beta}(u),\ \ \forall\ x,u.
\end{split}
\end{equation}
\end{Remark}

\begin{Remark}\label{Remark1b}
From \cite{Praly1}, a system of the form \eqref{system2a} that admits a $C^1$ function $V$ satisfying the inequalities (\ref{Vx1}) and (\ref{dotVx1}) is called exp-ISS, and the function
$V$ is called an exp-ISS Lyapunov function of \eqref{system2a}.  Moreover, by Proposition 8 of \cite{Praly1} or Theorem 3 of \cite{Sontag5},  system \eqref{system0a} is ISS if and only if it is exp-ISS. In this paper, we further call \eqref{system2a} strong exp-ISS
if it admits a $C^1$ function $\bar{V}$ satisfying
the inequalities (\ref{bVx1}) and (\ref{dotbVx3a}), and call $\bar{V}$ a strong exp-ISS Lyapunov function of \eqref{system2a}. Thus, Lemma \ref{Lemma0} shows that exp-ISS is equivalent to strong exp-ISS. It will be seen in the proof of Theorem \ref{Theorem1}
that Lemma \ref{Lemma0} plays the key role to eliminate the growth condition of the nonlinear functions $G_i (Z_i, e_i, d (t))$.
\end{Remark}

\section{Main Result}\label{Result}
In this section,  by combining Lemma \ref{Lemma0} with the multiple Lyapunov functions and average dwell time method, we will design a distributed switched output feedback control law to solve the cooperative global robust stabilization problem for system \eqref{system1a}.

For convenience, let $Z\!=\!\mbox{col}(Z_{1},\cdots,Z_{N})$, $e\!=\!\mbox{col}(e_{1},\cdots$, $e_{N})$, $e_{v\sigma(t)}\!=\!\mbox{col}(e_{v\sigma(t)1},\cdots,e_{v\sigma(t)N})$. From \eqref{evi}, we have $e_{v\sigma(t)} $ $= H_{\sigma(t)} e $ for all $ t \geq 0$. Then we can put the closed loop system composed of (\ref{system1a}) and (\ref{ui2a}) into the following form
\begin{equation}\label{xc1a}
\begin{split}
\dot{Z}&=F(Z,e,d(t))\\
\dot{e}&=\bar{G}_{\sigma(t)}(Z,e,d(t))\\
\end{split}
\end{equation}
where $F(Z,e,d(t))=\mbox{col}\big(F_{1}(Z_{1},e_{1},d(t)),\cdots,F_{N}(Z_{N}$, $e_{N},d(t))\big)$, $\bar{G}_{\sigma(t)}(Z,e,d(t))=\mbox{col}\big(G_{1}(Z_{1},e_{1},d(t))+b_{1}k_{1}(e_{v\sigma(t)1}),\cdots,G_{N}(Z_{N},e_{N},d(t))+b_{N}k_{N}(e_{v\sigma(t)N})\big)$.
 Let
 $f_{c\sigma(t)}(x_{c},d(t))=\mbox{col}\big(F(Z,e,d(t)), \bar{G}_{\sigma(t)}(Z,e,d(t))\big)$ and $x_{c}=\mbox{col}(Z,e)$. Then we can further put system \eqref{xc1a} into the following compact form
\begin{equation}\label{xc2a}
\begin{split}
\dot{x}_{c}=f_{c\sigma(t)}(x_{c},d(t)).
\end{split}
\end{equation}
Since (\ref{ui2a}) is a switched control law,  the closed loop system \eqref{xc2a} is a switched nonlinear system. To analyze the stability of system \eqref{xc2a}, we will resort to the multiple Lyapunov functions and average dwell time method \cite{Guo1,Liberzon1}.
From Theorem 4 of \cite{Guo1},   if the closed loop system
(\ref{xc2a}) admits   multiple $C^{1}$ Lyapunov functions
$U_{p}(x_{c})$, $p \in \mathcal{P}$, satisfying
  \begin{equation}\label{Up1a}
\begin{split}
\tilde{\alpha}_{1}(\|x_{c}\|)\leq
U_{p}(x_{c})\leq\tilde{\alpha}_{2}(\|x_{c}\|),\ \ \forall x_{c}, \ \
\forall p\in \mathcal{P}
\end{split}
\end{equation}
 \begin{equation}\label{dotUp1a}
\begin{split}
\dfrac{\partial U_{p}}{\partial
x_{c}}f_{cp}(x_{c},d(t))\leq-\lambda_{0}U_{p}(x_{c}),\ \ \forall
x_{c},\ \ \forall p\in \mathcal{P}
\end{split}
\end{equation}
for some class $\mathcal{K}_{\infty}$ functions
$\tilde{\alpha}_{1}(\cdot)$ and $\tilde{\alpha}_{2}(\cdot)$, and
some positive numbers $\lambda_{0}$, then the origin of
(\ref{xc2a})  is globally asymptotically stable for every $ \sigma \in S_{ave}[\tau_{d},N_{0}]$ with $\tau_{d}>\frac{\ln
\mu_{0}}{\lambda_{0}}$ and arbitrary $N_{0}$,  where $\mu_{0}=\sup\limits_{x_{c}\neq0}\frac{\tilde{\alpha}_{2}(\|x_{c}\|)}{\tilde{\alpha}_{1}(\|x_{c}\|)}$.



 \begin{Assumption}\label{Ass1.3}
 For a given compact subset $\mathbb{D} \subset \mathbb{R}^{n_{d}}$, the subsystem $\dot{Z}_{i}=F_{i}(Z_{i},e_{i},d(t))$ admits  a $C^{1}$  exp-ISS Lyapunov function $V_{i}(Z_{i})$ such that, for any $d(t)\in\mathbb{D}$,
 \begin{equation}\label{Vzi1a}
\begin{split}
\check{\alpha}_{1i}(\|Z_{i}\|)\leq
V_{i}(Z_{i})\leq\check{\alpha}_{2i}(\|Z_{i}\|),\ \
\forall\ Z_{i}
\end{split}
\end{equation}
 \begin{equation}\label{dotVzi1a}
\begin{split}
\dot{V}_{i}(Z_{i})\leq
-\lambda_{2}V_{i}(Z_{i})+\pi_{i}(e_{i}),\ \ \forall\
Z_{i},e_{i}
\end{split}
\end{equation}
 where  $\check{\alpha}_{1i}(\cdot)$ and $\check{\alpha}_{2i}(\cdot)$ are some class $\mathcal{K}_\infty$ functions with $\check{\alpha}_{1i}(\cdot)$ satisfying $\lim\limits_{s\rightarrow0^{+}}\sup\frac{s^{2}}{\check{\alpha}_{1i}(s)}<\infty$,  $\lambda_{2}$ is some  positive real number, and $\pi_{i}(\cdot)$ is some smooth positive definite function.
 \end{Assumption}
 \begin{Remark}\label{Remark1e}
Assumption \ref{Ass1.3} implies that the subsystem
$\dot{Z}_{i}=F_{i}(Z_{i},e_{i},d(t))$ is exponentially
input-to-state stable with $e_{i}$ as the input and
$V_{i}(Z_{i})$ is an exp-ISS Lyapunov function of the subsystem
$\dot{Z}_{i}=F_{i}(Z_{i},e_{i},d(t))$.
\end{Remark}

Now we describe our main result as follows.
 \begin{Theorem}\label{Theorem1}
 Under Assumption \ref{Ass1.3}, for every $ \sigma \in S_{ave}[\tau_{d},N_{0}]$ with
$\tau_{d}>\frac{\ln \mu_{0}}{\lambda_{0}}$  and  arbitrary $N_{0}$,
there exists a distributed switched output feedback control law of the form
\begin{equation}\label{ui3a}
\bar{u}_{i}=-\rho_{i}(e_{v\sigma(t)i})e_{v\sigma(t)i}, \ \ i=1,\cdots,N
\end{equation}
where $\rho_{i}(\cdot)$, $i=1,\cdots,N$, are some sufficiently
smooth positive functions, $\lambda_{0}$ and $\mu_{0}$ are some positive real numbers, that solves the cooperative global robust stabilization problem of system
(\ref{system1a}).
\end{Theorem}
\begin{Proof}
First note that, since $G_{i}(Z_{i},e_{i},d(t))$
is smooth and $G_{i}(0,0,d(t))=0$ for all $d(t)\in\mathbb{R}^{n_{d}}$, by
Lemma 7.8 of \cite{Huang1}, there exist some smooth positive
definite functions $\gamma_{i}(Z_{i})$, $\chi_{i}(e_{i})$, such
that, for all $Z_{i}\in\mathbb{R}^{n_{i}}$, $e_{i}\in \mathbb{R}$
and $d(t)\in\mathbb{D}$,
\begin{equation}\label{gi1a}
\begin{split}
  |G_{i}(Z_{i},e_{i},d(t))|^{2}\leq&\gamma_{i}(Z_{i})+\chi_{i}(e_{i}).
\end{split}
\end{equation}

Without loss of generality, assume $b_{m} \leq b_i \leq b_{M}$.
By Lemma \ref{Lemma0} and Remark \ref{Remark1}, under Assumption \ref{Ass1.3},  the subsystem $\dot{Z}_{i}=F_{i}(Z_{i},e_{i},d(t))$ admits a $C^{1}$  strong exp-ISS Lyapunov function
$\bar{V}_{i}(Z_{i})$  such that, for any $d (t) \in \mathbb{D}$,
   \begin{equation}\label{Vzi2a}
\begin{split}
\bar{\alpha}_{1i}(\|Z_{i}\|)\leq
\bar{V}_{i}(Z_{i})\leq\bar{\alpha}_{2i}(\|Z_{i}\|),\ \
\forall\ Z_{i}
\end{split}
\end{equation}
\begin{equation}\label{dotbVzi1a}
\begin{split}
\dot{\bar{V}}_{i}(Z_{i})\leq\!-\!\big(\lambda_{0}\bar{V}_{i}(Z_{i})\!+\!c_{0}\gamma_{i}(Z_{i})\big)\!+\!\bar{\pi}_{i}(e_{i}),~ \forall~ Z_{i}, e_i,
\end{split}
\end{equation}
for some class $\mathcal{K}_{\infty}$ functions $\bar{\alpha}_{1i}(\cdot)$ and $\bar{\alpha}_{2i}(\cdot)$,
 some  positive real numbers $\lambda_{0}$ and $c_{0}$, and some
positive definite smooth function $\bar{\pi}_{i}(\cdot)$.

For each $p\in \mathcal{P}$, let $e_{vp}=H_{p}e$ and $G(Z,e,d(t))=\mbox{col} (G_{1}(Z_{1},e_{1},d(t)),$ $\cdots,G_{N}(Z_{N},e_{N},d(t)))$.  Then, from \eqref{system1a}, we have
\begin{equation}\label{dotevp1a}
\begin{split}
\dot{e}_{vp}&=H_{p}\dot{e}=H_{p}G(Z,e,d(t))+H_{p}B\bar{u}\\
&=\check{G}_{p}(Z,e_{vp},d(t))+H_{p}B\bar{u}
\end{split}
\end{equation}
where $B=\mbox{diag}(b_{1},\cdots,b_{N})$, $\bar{u}=\mbox{col}(\bar{u}_{1},\cdots,\bar{u}_{N})$, and $\check{G}_{p}(Z,e_{vp},d(t))=\mbox{col}(\check{G}_{p1}(Z,e_{vp},d(t)),\cdots,\check{G}_{pN}(Z,e_{vp}$, $d(t)))=H_{p}G(Z,H_{p}^{-1}e_{vp},d(t))$.
Since $H_{p}=[h_{ij}^{p}]$ is a constant $\mathcal{M}$ matrix, then by Lemma 2.5.3 of \cite{Horn1}, there exists a positive definite diagonal matrix $D_{p}=\mbox{diag}(d_{p1},\cdots,d_{pN})$ such that $D_{p}H_{p}\!+\!H_{p}^{T}D_{p}$ is positive definite. Thus $B(D_{p}H_{p}$ $+H_{p}^{T}D_{p})B$ is also a positive definite matrix. Define $d_{M}=\max\limits_{p\in\mathcal{P}}\{\max\limits_{i=1,\cdots,N}d_{pi}\}$, $d_{m}=\min\limits_{p\in\mathcal{P}}\{\min\limits_{i=1,\cdots,N}d_{pi}\}$ and $\tilde{\lambda}_{1}=$ $\min\limits_{p\in\mathcal{P}}\{\lambda_{\min}(B(D_{p}H_{p}+$ $H_{p}^{T}D_{p})B)\}$.

Let $\rho_{i}(e_{vpi})=k\omega_{i}(e^{2}_{vpi})$, where $k$ is a positive real number and $\omega_{i}(\cdot)\geq1$, $i=1,2,\cdots,N$, are some smooth non-decreasing functions to be determined later. Let
\begin{equation}\label{V0a}
\begin{split}
V_{ep}(e)=\sum_{i=1}^{N}d_{pi}b_{i}\int_{0}^{e^{2}_{vpi}}\omega_{i}(s)ds.
\end{split}
\end{equation}
Then,  it can be seen that $V_{ep}(e)$ is positive definite and radially unbounded. Thus,  
there  exist some class $\mathcal{K}_{\infty}$ functions $\hat{\beta}_{1p}(\cdot)$ and $\hat{\beta}_{2p}(\cdot)$ such that
\begin{equation}\label{V2a}
\begin{split}
\hat{\beta}_{1p}(\|e\|)\leq V_{ep}(e) \leq\hat{\beta}_{2p}(\|e\|), \ \forall e, \ \forall p\in\mathcal{P}.
\end{split}
\end{equation}
Choose two class $\mathcal{K}_{\infty}$ functions $\hat{\beta}_{1}(\cdot)$ and $\hat{\beta}_{2}(\cdot)$ such that, for all $p\in\mathcal{P}$ and all $e\in\mathbb{R}^{N}$,  $\hat{\beta}_{1}(\|e\|)\leq\hat{\beta}_{1p}(\|e\|)$ and $\hat{\beta}_{2}(\|e\|)\geq\hat{\beta}_{2p}(\|e\|)$.  Then
\begin{equation}\label{V3a}
\begin{split}
\hat{\beta}_{1}(\|e\|)\leq V_{ep}(e) \leq\hat{\beta}_{2}(\|e\|), \ \forall e, \ \forall p\in\mathcal{P}.
\end{split}
\end{equation}
Let $e_{vp}^{*}=\mbox{col}(\omega_{1}(e^{2}_{vp1})e_{vp1},\cdots,\omega_{N}(e^{2}_{vpN})e_{vpN})$. Then, by  \eqref{ui3a}, \eqref{gi1a}, and \eqref{dotevp1a}, for any $\varepsilon_{2} >0$ and any $d(t)\in\mathbb{D}$, the time derivative of $V_{ep}(e)$ is given by
\begin{equation}\label{Vep1a}
\begin{split}
&\dfrac{\partial V_{ep}}{\partial e}(G(Z,e,d(t))+B\bar{u})\\
=&\dfrac{\partial V_{ep}}{\partial e_{vp}}(\check{G}_{p}(Z,e_{vp},d(t))+H_{p}B\bar{u})\\
 =&2\sum_{i=1}^{N}d_{pi}b_{i}\omega_{i}(e^{2}_{vpi})e_{vpi}\bigg(\check{G}_{pi}(Z,e_{vp},d(t))\\
 &-\sum_{j=1}^{N}h_{ij}^{p}b_{j}k\omega_{j}(e^{2}_{vpj})e_{vpj}\bigg)\\
  \leq&\varepsilon_{2}\sum_{i=1}^{N}d^{2}_{M}b^{2}_{M}\omega_{i}^{2}(e^{2}_{vpi})e^{2}_{vpi}+\dfrac{1}{\varepsilon_{2}}\sum_{i=1}^{N}|\check{G}_{pi}(Z,e_{vp},d(t))|^{2}\\
 &-k(e_{vp}^{*})^{T}B(D_{p}H_{p}+H_{p}^{T}D_{p})Be_{vp}^{*}\\
 \leq&-\sum_{i=1}^{N}(k\tilde{\lambda}_{1}-\varepsilon_{2}d^{2}_{M}b^{2}_{M})\omega_{i}^{2}(e^{2}_{vpi})e^{2}_{vpi}\\
 &+\dfrac{1}{\varepsilon_{2}}\|\check{G}_{p}(Z,e_{vp},d(t))\|^{2}\\
 \leq&-\sum_{i=1}^{N}(k\tilde{\lambda}_{1}-\varepsilon_{2}d^{2}_{M}b^{2}_{M})\omega_{i}^{2}(e^{2}_{vpi})e^{2}_{vpi}\\
 &+\dfrac{\|H_{p}\|^{2}}{\varepsilon_{2}}\sum_{i=1}^{N}|G_{i}(Z_{i},e_{i},d(t))|^{2}\\
 \leq&-\sum_{i=1}^{N}(k\tilde{\lambda}_{1}-\varepsilon_{2}d^{2}_{M}b^{2}_{M})\omega_{i}^{2}(e^{2}_{vpi})e^{2}_{vpi}\\
 &+\dfrac{\|H_{p}\|^{2}}{\varepsilon_{2}}\sum_{i=1}^{N}(\gamma_{i}(Z_{i})+\chi_{i}(e_{i}))\\
\end{split}
\end{equation}
Let $\bar{V}(Z)=\sum_{i=1}^{N}\bar{V}_{i}(Z_{i})$. By inequality \eqref{dotbVzi1a} we have
 \begin{equation}\label{dotVz0a}
\begin{split}
\dot{\bar{V}}(Z)\leq\!-\!\sum_{i=1}^{N}\big(\lambda_{0}\bar{V}_{i}(Z_{i})\!+\!
c_{0}\gamma_{i}(Z_{i})\big)\!+\!\sum_{i=1}^{N}\bar{\pi}_{i}(e_{i}).
 \end{split}
\end{equation}
Finally, let $U_{p}(x_{c})=\bar{V}(Z)+V_{ep}(e)$. 
Clearly, there exist two class $\mathcal{K}_{\infty}$ functions $\tilde{\alpha}_{1}(\cdot)$ and $\tilde{\alpha}_{2}(\cdot)$ such that
the condition (\ref{Up1a}) is satisfied for all $p\in \mathcal{P}$. Also, according to \eqref{Vep1a} and \eqref{dotVz0a}, we have
\begin{equation}\label{dotUp2a}
\begin{split}
 &\dfrac{\partial U_{p}}{\partial x_{c}}f_{cp}(x_{c},d(t))\\
  \leq & -\sum_{i=1}^{N}\lambda_{0}\bar{V}_{i}(Z_{i})- \sum_{i=1}^{N}(c_{0}-\dfrac{\|H_{p}\|^{2}}{\varepsilon_{2}})\gamma_{i}(Z_{i})\\
  &-\sum_{i=1}^{N}(k\tilde{\lambda}_{1}-\varepsilon_{2}d^{2}_{M}b^{2}_{M})\omega_{i}^{2}(e^{2}_{vpi})e^{2}_{vpi}+\hat{\rho}(e)\\
\end{split}
\end{equation}
with
$\hat{\rho}(e)=\sum_{i=1}^{N}\big(\frac{\|H_{p}\|^{2}}{\varepsilon_{2}}\chi_{i}(e_{i})+\bar{\pi}_{i}(e_{i})\big)$.
By lemma $7.8$ of \cite{Huang1} again, for each $p\in \mathcal{P}$,
there exist some smooth positive functions
$\tilde{\rho}_{ip}(e_{vpi})\geq1$, $i=1,\cdots,N$ such that
$\hat{\rho}(e)=\hat{\rho}(H_{p}^{-1}e_{vp})\leq\sum_{i=1}^{N}\tilde{\rho}_{ip}(e_{vpi})e_{vpi}^{2}$.
Choose some smooth functions $\tilde{\rho}_{i}(e_{vpi})$ such that $\tilde{\rho}_{i}(e_{vpi})\geq\tilde{\rho}_{ip}(e_{vpi})$ for all
$p\in \mathcal{P}$. Then $\hat{\rho}(e)\leq\sum_{i=1}^{N}\tilde{\rho}_{i}(e_{vpi})e_{vpi}^{2}$ for all
$p\in \mathcal{P}$.
Define a function $\Delta_{i}: \mathbb{R}^{+}\rightarrow \mathbb{R}$ such that $\Delta_{i}(s)=\max\limits_{0\leq|e_{vpi}|\leq s}\tilde{\rho}_{i}(e_{vpi})$. Clearly, $\Delta_{i}(\cdot)$ is non-decreasing and $\Delta_{i}(|e_{vpi}|)\geq\tilde{\rho}_{i}(e_{vpi})\geq1$ for any $e_{vpi}\in \mathbb{R}$.
Then we choose the positive real number $k$ such that
\begin{equation}\label{k1a}
\begin{split}
k\geq \dfrac{1}{\tilde{\lambda}_{1}}(\varepsilon_{2}d^{2}_{M}b^{2}_{M}+\lambda_{0}d_{M}b_{M}+1)
\end{split}
\end{equation}
and choose $\varepsilon_{2}$ such that $c_{0}-\frac{\|H_{p}\|^{2}}{\varepsilon_{2}}\geq0$. Then,
\begin{equation}\label{dotUp3a}
\begin{split}
  &\dfrac{\partial U_{p}}{\partial x_{c}}f_{cp}(x_{c},d(t))\\
   \leq& -\lambda_{0}\sum_{i=1}^{N}\bar{V}_{i}(Z_{i}) -\lambda_{0}\sum_{i=1}^{N}d_{M}b_{M}\omega_{i}^{2}(e^{2}_{vpi})e^{2}_{vpi}\\
   &-\sum_{i=1}^{N}(\omega_{i}^{2}(e^{2}_{vpi})-\tilde{\rho}_{i}(e_{vpi}))e^{2}_{vpi},\ \ \forall p\in \mathcal{P}
\end{split}
\end{equation}
Finally, choose the smooth non-decreasing function $\omega_{i}(\cdot)$ such that
\begin{equation}\label{omega1a}
\begin{split}
\omega_{i}(e_{vpi}^{2})\geq \Delta_{i}\bigg(\dfrac{1+e_{vpi}^{2}}{2}\bigg).
\end{split}
\end{equation}
Then, since $\Delta_{i}\big(\frac{1+e_{vpi}^{2}}{2}\big)\geq\Delta_{i}(|e_{vpi}|)\geq\tilde{\rho}_{i}(e_{vpi})\geq1$, we have
\begin{equation}\label{omega2a}
\begin{split}
\omega_{i}^{2}(e_{vpi}^{2})\geq\omega_{i}(e_{vpi}^{2})\geq \Delta_{i}\bigg(\dfrac{1+e_{vpi}^{2}}{2}\bigg)\geq\tilde{\rho}_{i}(e_{vpi})
\end{split}
\end{equation}
and thus
\begin{equation}\label{dotUp4a}
\begin{split}
  &\dfrac{\partial U_{p}}{\partial x_{c}}f_{cp}(x_{c},d(t))\\
   \leq& -\lambda_{0}\sum_{i=1}^{N}\bar{V}_{i}(Z_{i}) -\lambda_{0}\sum_{i=1}^{N}d_{M}b_{M}\omega_{i}(e^{2}_{vpi})e^{2}_{vpi}\\
\leq& -\lambda_{0}\bar{V}(Z)-\lambda_{0}V_{ep}(e)  =-\lambda_{0}U_{p},\ \forall p\in \mathcal{P}.\\
\end{split}
\end{equation}
Let $\mu_{0}=\sup\limits_{x_{c}\neq0}\frac{\tilde{\alpha}_{2}(\|x_{c}\|)}{\tilde{\alpha}_{1}(\|x_{c}\|)}$.
The proof is thus completed by invoking Theorem 4 of \cite{Guo1}
as rephrased at the beginning of this section.
\end{Proof}
\begin{Remark}\label{Remark3.2}
From the proof of Theorem \ref{Theorem1}, we know that if we replace the exp-ISS condition \eqref{dotVzi1a} in Assumption \ref{Ass1.3} by the strong exp-ISS condition \eqref{dotbVzi1a}, then the condition $\lim\limits_{s\rightarrow0^{+}}\sup\frac{s^{2}}{\check{\alpha}_{1i}(s)}<\infty$ in Assumption \ref{Ass1.3} can be removed.
\end{Remark}

\begin{Remark}
The recent paper \cite{Su6}  handled the cooperative output regulation problem for nonlinear systems with static directed graph. The Lyapunov function (\ref{V0a}) here is similar to
that used in Lemma 5 of \cite{Su6}.
\end{Remark}
\begin{Example}
Consider the following controlled Lorenz multi-agent
systems taken from \cite{Dongyi1}:
\begin{equation}\label{Lorenzsystem2}
\begin{split}
 &\dot{z}_{1i}=-L_{1i}z_{1i}+L_{1i}e_{i}\\
 &\dot{z}_{2i}=L_{2i}z_{2i}+z_{1i}e_{i}\\
 &\dot{e}_{i}= L_{3i}z_{1i}-e_{i}-z_{1i}z_{2i} + b_{i}\bar{u}_{i},~i=1,2,3
 \end{split}
\end{equation}
where $b_{i}=1$, $L_{i}=\mbox{col} (L_{1i}, L_{2i}, L_{3i})$ is a constant
parameter vector that satisfies $L_{1i}>0$, $L_{2i}<0$ and
$L_{3i}>0$. Clearly, system \eqref{Lorenzsystem2} is in the form of \eqref{system1a} with $Z_{i}=\mbox{col}(z_{1i},z_{2i})$.
To account for the uncertainty, for $i=1,2,3$, let $L_{i}=\bar{L}_{i}+d_{i}$ where
$\bar{L}_{i}=\mbox{col}(\bar{L}_{1i},\bar{L}_{2i},\bar{L}_{3i})=\mbox{col}(3,-3.2,1.6)$
denotes the nominal value of $L_{i}$, and $d_{i}=\mbox{col}(d_{1i},d_{2i},d_{3i})$
represents the uncertainty of $L_{i}$. Assume $d=\mbox{col}(d_{1},d_{2},d_{3})\in\mathbb{D}=\{d\in\mathbb{R}^{9}|~|d_{ji}|\leq0.2,i,j=1,2,3\}$. Let $H_{1}=\left[
                                     \begin{array}{ccc}
                                       2 & 0 & -1 \\
                                       -1 & 1 & 0 \\
                                       -1 & 0 & 2 \\
                                     \end{array}
                                   \right]
$, $H_{2}=\left[
                                     \begin{array}{ccc}
                                       1 & 0 & 0 \\
                                       0 & 2 & -1 \\
                                       -1 & -1 & 2 \\
                                     \end{array}
                                   \right]
$ which are both $\mathcal{M}$ matrices. Next, let a Lyapunov function candidate be $V_{i}(Z_{i})=\frac{1}{2}z_{1i}^{2}+\frac{1}{4}z_{1i}^{4}+\frac{1}{2}z_{2i}^{2}$. Then it is possible to show that $\dot{V}_{i}(Z_{i})\leq-2.3z_{1i}^{2}-1.8z_{1i}^{4}-2.5z_{2i}^{2}+5.2e_{i}^{2}+26.5e_{i}^{4}\leq-4.6V_{i}(Z_{i})+5.2e_{i}^{2}+26.5e_{i}^{4}$. Thus Assumption \ref{Ass1.3} is satisfied.


\vspace{-0.2cm}
\begin{figure}[H]
\centering
\includegraphics[scale=0.5]{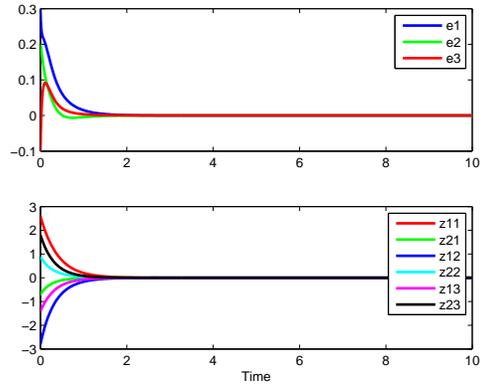}
\caption{State trajectories} \label{states}
\end{figure}
\vspace{-0.2cm}

As a result, using the procedure introduced in this section,
the cooperative robust stabilization problem for
system \eqref{Lorenzsystem2} is solvable by a switched control law:
$\bar{u}_{i}=-12(e_{v\sigma(t)i}^{4}+1)e_{v\sigma(t)i},~i=1,2,3$, for any $\sigma(t)\in S_{ave}[\tau_{d},N_{0}]$ with $\tau_{d}>2.72$ and arbitrary $N_{0}$.
Simulation is conducted for the following  specific switching signal
\begin{equation}\label{sigmat}
\begin{split}
\sigma(t)= \left\{
  \begin{array}{ll}
    1, & \hbox{if\ $sT\leq t< (s+\dfrac{1}{2})T$} \\
    2, & \hbox{if\ $(s+\dfrac{1}{2})T\leq t< (s+1)T$} \\
  \end{array}
\right.
 \end{split}
\end{equation}
where $s=0,1,2,\cdots$, and $T=2\tau_{d} = 6$ sec. Other data for the simulation are $d_{i}=[-0.2,0.1,-0.2]^{T}$ for $i=1,2,3$, $Z_{1}(0)=[2.6, -0.7]^{T}$, $Z_{2}(0)=[-2.8, 0.9]^{T}$, $Z_{3}(0)=[-1.4, 1.8]^{T}$ and $e(0)=[0.3, 0.2, -0.1]^{T}$.
Figure \ref{states} shows the state trajectories of the closed-loop system.
Due to the space limit, the details for designing the controller are omitted.

\end{Example}

\section{An Application}\label{Application}
In the past few years, the cooperative control problems for nonlinear multi-agent systems have been extensively studied for the static network case in \cite{Ding1,Dongyi1,Isidori1,Su4,Su6,XuHong1}, and for the switching network case in \cite{Liu2,Zhao1}. Note that the nonlinear systems considered in \cite{Zhao1} contain no disturbance and uncertainty,  and the nonlinear systems considered in \cite{Liu2} need to satisfy certain growth condition. Also, the switching network is assumed to be undirected in \cite{Liu2}.
In this section, we will apply our main result to the  cooperative global robust output regulation problem for
a class of nonlinear multi-agent systems with switching network. The systems studied here contain both disturbance and parameter uncertainty, and do not need to satisfy the growth condition, moreover, the switching network is relaxed to be directed.

 Consider the nonlinear multi-agent systems in normal form with unity relative degree as follows, which is the same as those studied in \cite{Dongyi1,Liu2}.
\begin{equation}\label{system1}
\begin{split}
  \dot{z}_{i} =& f_{i}(z_{i},y_{i},v,w)\\
  \dot{y}_{i} =& b_{i}(w)u_{i}+g_{i}(z_{i},y_{i},v,w)\\
  e_{i}=&y_{i}-q(v,w),\ \ \ i=1,\cdots,N
\end{split}
\end{equation}
where, for $i=1,\cdots,N$, $(z_{i},y_{i})\in \mathbb{R}^{n_{z_{i}}}\times\mathbb{R}$ is the state, $u_{i}\in \mathbb{R}$ is the input, $e_{i}\in \mathbb{R}$ is the error output,
$w\in\mathbb{R}^{n_{w}}$ is an uncertain parameter vector, and $v(t)\in\mathbb{R}^{n_{v}}$ is an exogenous signal representing both reference input and disturbance.
It is assumed that $v(t)$ is generated by a linear system of the following form
\begin{equation}\label{exosystem1}
\begin{split}
  \dot{v} = Sv,\ \ y_{0}=q(v,w)
\end{split}
\end{equation}
and all functions in (\ref{system1}) and (\ref{exosystem1}) are
globally defined, sufficiently smooth, and satisfy
$f_{i}(0,0,0,w)=0$, $g_{i}(0,0,0,w)=0$ and $q(0,w)=0$ for all
$w\in \mathbb{R}^{n_{w}}$.


As in \cite{Liu2}, the plant (\ref{system1}) and the exosystem (\ref{exosystem1}) together can be viewed as a multi-agent system of $(N
+1)$ agents with (\ref{exosystem1}) as the leader and the $N$
subsystems of (\ref{system1}) as $N$ followers. Given the plant (\ref{system1}),  the exosystem (\ref{exosystem1}),
and a switching signal $\sigma(t)$,  we can define a time-varying
digraph
$\bar{\mathcal{G}}_{\sigma(t)}=(\bar{\mathcal{V}},\bar{\mathcal{E}}_{\sigma(t)})$,
where $\bar{\mathcal{V}}=\{0,1\cdots,N\}$ with $0$ associated with
the leader system and with $ i = 1,\cdots,N$ associated with the $N$
followers, respectively, and $\bar{\mathcal{E}}_{\sigma(t)}\subseteq
\bar{\mathcal{V}}\times\bar{\mathcal{V}}$ for all $t\geq0$.
 For all
$t \geq 0$, each $j=0,1,\cdots,N$, $i=1,\cdots,N$, $i \neq j$,
$(j,i)\in\bar{\mathcal{E}}_{\sigma(t)} $  if and only if the control
$u_{i} (t)$ can make use of $y_{i} (t) -y_{j} (t)$ for feedback
control.
Let
$\mathcal{G}_{\sigma(t)}=(\mathcal{V},\mathcal{E}_{\sigma(t)})$ where $\mathcal{V}=\{1\cdots,N\}$,
$\mathcal{E}_{\sigma(t)}
\subseteq\mathcal{\mathcal{V}\times\mathcal{V}}$ is obtained from
$\bar{\mathcal{E}}_{\sigma(t)}$ by removing all edges between the
node $0$ and the nodes in $\mathcal{V}$ for all $t \geq 0$.

Let  $\bar{\mathcal{A}} (t)=[\bar{a}_{ij} (t)]\in \mathbb{R}^{(N+1)\times
(N+1)}$  be the adjacency matrix of the digraph $\bar{\mathcal{G}}_{\sigma (t)}$,
where $\bar{a}_{ii} (t)=0$ and $\bar{a}_{ij} (t)=1\Leftrightarrow
(j,i)\in\mathcal{\bar{E}}_{\sigma (t)}$, $i,j=0,1,\cdots,N$.
Define the virtual regulated output as follows:
\begin{equation}\label{evi1}
\begin{split}
e_{v\sigma(t)i}=\sum_{j=0}^{N}\bar{a}_{ij}(t)(y_{i}-y_{j}).
\end{split}
\end{equation}
Let
$e_{v\sigma(t)}=\mbox{col} (e_{v\sigma(t)1},\cdots,e_{v\sigma(t)N})$, $ e = \mbox{col}
(e_{1},\cdots,e_{N})$,  $e_{0}=0$ and
$H_{ \sigma(t)}=[h_{ij}(t)]_{i,j=1}^{N}$ with $h_{ii}(t) =\sum_{j=0}^{N}\bar{a}_{ij}  (t)$ and $h_{ij}(t)=-\bar{a}_{ij}  (t)$ for $i\neq j$.
It can be verified that $e_{v\sigma(t)}=H_{\sigma(t)}e$.
Then our control law will be of the following form
\begin{equation}\label{ui2}
\begin{split}
  &u_{i}=\hat{k}_{i}(\eta_{i},e_{v\sigma(t)i} ),~\dot{\eta}_{i}= \hat{g}_{i}(\eta_{i},e_{v\sigma(t)i})
\end{split}
\end{equation}
where the functions $\hat{k}_i$ and $\hat{g}_i$ are sufficiently smooth
vanishing at the origin.  A control law of the form (\ref{ui2}) is
called a distributed switched output feedback control law since $e_{v\sigma(t)i}$ is a switching signal and depends on $(y_i - y_j)$ if only if the node $j$ is a neighbor of the node $i$.
Then we describe our problem as follows:
%
\begin{Problem}\label{CORPS}
 Given the multi-agent system (\ref{system1}), the exosystem (\ref{exosystem1}),
a group of  digraphs $\bar{\mathcal{G}}_{p} = (\bar{\mathcal{V}},\bar{\mathcal{E}}_{p})$ with $\bar{\mathcal{V}}=\{0,1\cdots,N\}$, $p = 1,\cdots, n_{0}$,  and some compact subsets
$\mathbb{V}\subset\mathbb{R}^{n_{v}}$ and
$\mathbb{W}\subset\mathbb{R}^{n_{w}}$ with $0\in\mathbb{W}$ and
$0\in\mathbb{V}$, find $\tau_{d},N_{0}$, and  a control law of the
form (\ref{ui2})  such that, for any $v(t)\in \mathbb{V}$, $w\in
\mathbb{W}$, and any $ \sigma \in S_{ave}[\tau_{d},N_{0}]$, the
trajectory of the closed-loop system composed of \eqref{system1} and \eqref{ui2} starting from any
initial state $z_{i}(0)$, $y_{i}(0)$ and $\eta_{i}(0)$ exists and is bounded for all $t\geq0$, and
$\lim_{t \to \infty}e(t)=0$.
\end{Problem}


The above problem was studied recently in \cite{Liu2} where it was shown that this problem can be converted to the problem studied in Section \ref{Result} under the following assumptions.

\begin{Assumption}\label{Ass2.1}
The exosystem is neutrally stable, i.e., all the eigenvalues of $S$ are semi-simple with zero real parts.
\end{Assumption}

\begin{Assumption}\label{Ass2.2}
For $i=1,\cdots,N$,  $|b_{i}(w)|>0$ for all $w\in\mathbb{R}^{n_{w}}$.
\end{Assumption}


\begin{Assumption}\label{Ass2.3}
There exist globally defined smooth functions $\mathbf{z}_{i}:\mathbb{R}^{n_{v}}\times\mathbb{R}^{n_{w}}{\rightarrow}\mathbb{R}^{n_{z_{i}}}$ with $\mathbf{z}_{i}(0,w)=0$ such that
\begin{equation}\label{Zi}
\begin{split}
  \dfrac{\partial\mathbf{z}_{i}(v,w)}{\partial v}Sv=f_{i}(\mathbf{z}_{i}(v,w),q(v,w),v,w)
\end{split}
\end{equation}
for all $(v,w)\in\mathbb{R}^{n_{v}}\times\mathbb{R}^{n_{w}}$, $i=1,\cdots,N$.
\end{Assumption}

\begin{Assumption}\label{Ass2.4}
 Let $\mathbf{u}_{i}(v,w)=b_{i}^{-1}(\dfrac{\partial q(v,w)}{\partial v}Sv
  -g_{i}(\mathbf{z}_{i}(v,w),q(v,w),v,w))$, $i=1,\cdots,N$. Then
  $\mathbf{u}_{i}(v,w)$  are polynomials in $v$ with coefficients depending on $w$.
\end{Assumption}

Under Assumptions \ref{Ass2.1} to \ref{Ass2.4},  there exist integers $s_i$, $i = 1, \cdots, N$, such that,  for
any Hurwitz matrices $M_{i}\in \mathbb{R}^{s_{i}\times s_{i}}$,  and  any column vector $N_{i}\in \mathbb{R}^{s_{i}\times 1}$ with
$(M_{i},N_{i})$ controllable, the following linear dynamic compensator
\begin{equation}\label{doteta1}
\begin{split}
  \dot{\eta}_{i}=M_{i}\eta_{i}+N_{i}u_{i},\  i = 1, \cdots, N,
\end{split}
\end{equation}
 is an internal model of system \eqref{system1} \cite{Liu2}. Moreover, there exist some functions  $\theta_i:  {R}^{n_v+n_w}{\rightarrow} {R}^{s_i}$ vanishing at the origin, and some row vectors $\Psi_{i}\in \mathbb{R}^{1 \times s_{i}}$ such that  the following
coordinate and input transformation on the internal model (\ref{doteta1}) and the plant (\ref{system1})
\begin{equation}\label{transformation}
\begin{split}
  \bar{z}_{i}=&z_{i}-\mathbf{z}_{i}(v,w),\ \ \tilde{\eta}_{i}=\eta_{i}-\theta_{i}(v,w)-N_{i}b_{i}^{-1}e_{i}\\
  e_{i}=&y_{i}-q(v,w),\ \ \bar{u}_{i}=u_{i}-\Psi_{i}\eta_{i},\ \ i=1,\cdots,N
\end{split}
\end{equation}
gives rise to the so-called augmented  system of the plant (\ref{system1}) and the exosystem (\ref{exosystem1}) as follows.
\begin{equation}\label{system2}
\begin{split}
  \dot{\bar{z}}_{i}=&\bar{f}_{i}(\bar{z}_{i},e_{i},d(t))\\
  \dot{\tilde{\eta}}_{i}=&M_{i}\tilde{\eta}_{i}+M_{i}N_{i}b_{i}^{-1}e_{i}-N_{i}b_{i}^{-1}\bar{g}_{i}(\bar{z}_{i},e_{i},d(t))\\
  \dot{e}_{i}=&\bar{g}_{i}(\bar{z}_{i},e_{i},d(t))+b_{i}\Psi_{i}\tilde{\eta}_{i}+\Psi_{i}N_{i}e_{i}+b_{i}\bar{u}_{i}
\end{split}
\end{equation}
where $d(t)=(v,w)$, $\bar{f}_{i}(\bar{z}_{i},e_{i},d(t))=f_{i}(\bar{z}_{i}+\mathbf{z}_{i},e_{i}+q,v,w)-f_{i}(\mathbf{z}_{i},q,v,w)$ and
$\bar{g}_{i}(\bar{z}_{i},e_{i},d(t))=g_{i}(\bar{z}_{i}+\mathbf{z}_{i},e_{i}+q,v,w)-g_{i}(\mathbf{z}_{i},q,v,w)$.
Clearly, $\bar{f}_{i}(0,0,d(t))=0$ and $\bar{g}_{i}(0,0,d(t))=0$ for any $d(t)\in\mathbb{R}^{n_{v}}\times\mathbb{R}^{n_{w}}$.

By the internal model principle as can be found from \cite{Liu2} or \cite{Dongyi1}, if a control law of the following form
\begin{equation}\label{barui2}
\begin{split}
\bar{u}_{i}=k_{i}(e_{v\sigma(t)i}),\ \ i=1,\cdots,N
\end{split}
\end{equation}
 solves the cooperative global robust stabilization problem of (\ref{system2}),
 then the cooperative global robust output regulation of system (\ref{system1}) is solved by the following distributed switched output feedback controller:
 \begin{equation}\label{ui4}
\begin{split}
u_{i}&=k_{i}({e}_{v\sigma(t)i})+\Psi_{i}\eta_{i}\\
\dot{\eta}_{i}&=M_{i}\eta_{i}+N_{i}u_{i},\ \ i=1,\cdots,N.
\end{split}
\end{equation}

It was further shown in \cite{Liu2} that the cooperative global robust stabilization problem of (\ref{system2}) was solvable by a control law of the form (\ref{barui2}) under the following three assumptions:

 \begin{Assumption}\label{Ass3.1o}
 For the given compact subset $\mathbb{D} \subset \mathbb{R}^{n_{v}}\times\mathbb{R}^{n_{w}}$, the subsystem $\dot{\bar{z}}_{i}=\bar{f}_{i}(\bar{z}_{i},e_{i},d(t)$ admits a $C^{1}$ function $V_{\bar{z}_{i}}(\bar{z}_{i})$ such that,  for any $d(t)\in\mathbb{D}$,
    \begin{equation}\label{Vzi}
\begin{split}
\alpha_{1i}(\|\bar{z}_{i}\|)\leq
V_{\bar{z}_{i}}(\bar{z}_{i})\leq\alpha_{2i}(\|\bar{z}_{i}\|),\ \
\forall\ \bar{z}_{i}
\end{split}
\end{equation}
\begin{equation}\label{dotVzi}
\begin{split}
\dot{V}_{\bar{z}_{i}}(\bar{z}_{i})\leq
-\lambda_{1} V_{\bar{z}_{i}}(\bar{z}_{i})+\beta_{i}(e_{i}),\ \ \forall\
\bar{z}_{i},e_{i}
\end{split}
\end{equation}
for some class $\mathcal{K}_{\infty}$ functions $\alpha_{1i}(\cdot)$ and $\alpha_{2i}(\cdot)$ with  $\alpha_{1i}(\cdot)$ satisfying $\lim\limits_{s\rightarrow0^{+}}\sup\frac{s^{2}}{{\alpha}_{1i}(s)}<\infty$, some  positive real number $\lambda_{1}$, and some  smooth  positive
definite function $\beta_{i}(\cdot)$.
 \end{Assumption}


\begin{Assumption}\label{Ass3.2o}
 For any $p\in\mathcal{P}$, $\mathcal{G}_{p}$ is undirected, and the node $0$ can reach every other node of the digraph $\bar{\mathcal{G}}_{p}$.
 \end{Assumption}

To introduce the last assumption, note that, since $\bar{g}_{i}(\bar{z}_{i},e_{i},d(t))$
is smooth and $\bar{g}_{i}(0,0,d(t))=0$ for all $d(t)\in\mathbb{D}$, by
Lemma 7.8 of \cite{Huang1}, there exist some smooth positive
definite functions $\delta_{i}(\bar{z}_{i})$, $l_{i}(e_{i})$ such
that, for all $\bar{z}_{i}\in\mathbb{R}^{n}$, $e_{i}\in \mathbb{R}$
and $d(t)\in\mathbb{D}$,
\begin{equation}\label{bargi}
\begin{split}
  |\bar{g}_{i}(\bar{z}_{i},e_{i},d(t))|^{2}\leq&\delta_{i}(\bar{z}_{i})+l_{i}(e_{i}).
\end{split}
\end{equation}

\begin{Assumption}\label{Ass3.3o} For some real number $c >0$, $c\delta_{i}(\bar{z}_{i}) \leq V_{\bar{z}_{i}}(\bar{z}_{i})$.
\end{Assumption}

Assumption \ref{Ass3.3o} is called a growth condition on the nonlinear function $\bar{g}_{i}(\bar{z}_{i},e_{i},d(t))$ which is quite restrictive, and Assumption \ref{Ass3.2o} requires
the graph $\mathcal{G}_{\sigma (t)}$ to be undirected for all $t \geq 0$ which may also be restrictive.
By making use of Theorem \ref{Theorem1} of this paper,  it is possible to remove Assumption \ref{Ass3.3o} and significantly relax \ref{Ass3.2o}.
For this purpose,  let $Z_{i}=\mbox{col} (\bar{z}_{i},\tilde{\eta}_{i})$, $G_{i}(Z_{i},e_{i},d(t))  = \bar{g}_{i}(\bar{z}_{i},e_{i},d(t))+b_{i}\Psi_{i}\tilde{\eta}_{i}+\Psi_{i}N_{i}e_{i}$ and $F_{i}(Z_{i},e_{i},d(t))=\mbox{col} (\bar{f}_{i}(\bar{z}_{i},e_{i},d(t)),M_{i}\tilde{\eta}_{i}+M_{i}N_{i}b_{i}^{-1}e_{i}-N_{i}b_{i}^{-1}\bar{g}_{i}(\bar{z}_{i},e_{i},d(t))$. Then the system (\ref{system2}) can be put in exactly the same form as \eqref{system1a}.  By Theorem 3.1, Assumption \ref{Ass3.2o} can be relaxed to the following

\begin{Assumption}\label{Ass3.2}
 For any $p\in\mathcal{P}$, every node $i=1,\cdots,N$ of the digraph $\bar{\mathcal{G}}_{p}$ is reachable from the node $0$.
 \end{Assumption}
\begin{Remark}
Under Assumption \ref{Ass3.2}, by lemma 4 in \cite{Hu1}, $H_{p}$ is an $\mathcal{M}$ matrix for any $p\in\mathcal{P}$.
\end{Remark}

We now further show that, by making use of Lemma \ref{Lemma0} and Remark \ref{Remark1}, Assumption \ref{Ass3.3o} can be removed. For this purpose, it suffices to show the following lemma.

 \begin{Lemma}\label{Lemma3.2}
 Under Assumption \ref{Ass3.1o}, 
 the subsystem $\dot{Z}_{i}=F_{i}(Z_{i},e_{i},d(t))$  admits a strong exp-ISS Lyapunov function $\bar{V}_{i}(Z_{i})$ such that,
 for any $d(t)\in\mathbb{D}$,
 \begin{equation}\label{VZi1}
\begin{split}
\bar{\alpha}_{1i}(\|Z_{i}\|)\leq
\bar{V}_{i}(Z_{i})\leq\bar{\alpha}_{2i}(\|Z_{i}\|),\ \ \forall\ Z_{i}
\end{split}
\end{equation}
\begin{equation} \label{dotVZi1}
\begin{split}
\dot{\bar{V}}_{i}(Z_{i})\leq&\!-\!(\lambda_{0}\bar{V}_{i}(Z_{i})\!+\!c_{0}\gamma_{i}(Z_{i}))+\bar{\pi}_{i}(e_{i}), \ \forall\ Z_{i}, e_i
\end{split}
\end{equation}
for some class $\mathcal{K}_{\infty}$ functions $\bar{\alpha}_{1i}(\cdot)$ and $\bar{\alpha}_{2i}(\cdot)$, some positive real numbers $\lambda_{0}$ and $c_{0}$, and some
positive definite smooth function $\bar{\pi}_{i}(\cdot)$.
 \end{Lemma}

\begin{Proof}
By Lemma \ref{Lemma0} and Remark \ref{Remark1}, under Assumption \ref{Ass3.1o},  
the subsystem $\dot{\bar{z}}_{i}=\bar{f}_{i}(\bar{z}_{i},e_{i},d(t))$ admits a  $C^{1}$ strong exp-ISS function $\bar{V}_{\bar{z}_{i}}(\bar{z}_{i})$ such that,
for any $d(t)\in\mathbb{D}$,
\begin{equation}\label{Vzi1}
\begin{split}
\hat{\alpha}_{1i}(\|\bar{z}_{i}\|)\leq
\bar{V}_{\bar{z}_{i}}(\bar{z}_{i})\leq\hat{\alpha}_{2i}(\|\bar{z}_{i}\|),\ \
\forall\ \bar{z}_{i}
\end{split}
\end{equation}
 \begin{equation}\label{dotVzi1}
\begin{split}
\dot{\bar{V}}_{\bar{z}_{i}}(\bar{z}_{i})\!\leq\!
-\!\big(\bar{\lambda}_{1}\bar{V}_{\bar{z}_{i}}(\bar{z}_{i})\!+\!
\bar{c}_{1}\delta_{i}(\bar{z}_{i})\big)\!+\!\bar{\beta}_{i}(e_{i}),\ \ \forall\
\bar{z}_{i}, e_i
\end{split}
\end{equation}
for some class $\mathcal{K}_{\infty}$ functions $\hat{\alpha}_{1i}(\cdot)$ and $\hat{\alpha}_{2i}(\cdot)$, some  positive real numbers $\bar{\lambda}_{1}$ and $\bar{c}_1$, and some
positive definite smooth function $\bar{\beta}_{i}(\cdot)$.

Let $V_{\tilde{\eta}_{i}}=\tilde{\eta}_{i}^{T}P_{i}\tilde{\eta}_{i}$, $\bar{V}_{i}(Z_{i})=\bar{V}_{\bar{z}_{i}}+V_{\tilde{\eta}_{i}}$,
 and $\gamma_{i}(Z_{i})=\delta_{i}(\bar{z}_{i})+\|\tilde{\eta}_{i}\|^{2}$,
 where $P_{i}$ is a symmetric positive definite matrix.
Then, by Lemma 3.1 and Remark 3.3 of \cite{Liu2},  $\bar{V}_{i}(Z_{i})$ satisfies both (\ref{VZi1}) and (\ref{dotVZi1}).
\end{Proof}

It is noted that conditions (\ref{VZi1}) and (\ref{dotVZi1}) are the same as conditions \eqref{Vzi2a} and \eqref{dotbVzi1a}) of Theorem \ref{Theorem1}.
Thus, combining  Theorem \ref{Theorem1}, Remark \ref{Remark3.2}, and Lemma \ref{Lemma3.2}  gives the following result.

\begin{Theorem}\label{Theorem3.1}
 Under Assumptions \ref{Ass2.1}-\ref{Ass3.1o}, and \ref{Ass3.2},  for every
 $ \sigma \in S_{ave}[\tau_{d}$, $N_{0}]$ with
$\tau_{d}>\frac{\ln \mu_{0}}{\lambda_{0}}$  and  arbitrary $N_{0}$, the cooperative global robust output regulation problem of system
(\ref{system1}) with the directed switching graph $\bar{ \mathcal{G}}_{\sigma(t)}$
is solvable by the distributed switched  output feedback control law of the
form (\ref{ui4}).
\end{Theorem}

\section{Conclusion}
In this paper, we have established a specific changing supply pair technique to analyze the exponential input to state stability for nonlinear systems. Then, combining this technique with multiple Lyapunov functions and the average dwell time method,  we have solved the cooperative global robust stabilization problem for a class of nonlinear multi-agent systems by a distributed switched output feedback control law. Finally, we have applied this result to the cooperative global robust output regulation problem for nonlinear multi-agent systems in normal form with unity relative degree under directed switching network.


\begin{thebibliography}{0}


\bibitem{Ding1}
Z. Ding (2013). Consensus output regulation of a class of heterogeneous nonlinear systems. \newblock {\em IEEE Transactions on Automatic Control}, 58(10), 2648-2653.

\bibitem{Dongyi1}
 Y. Dong and  J. Huang (2014). Cooperative global robust output regulation
for nonlinear multi-agent systems in output feedback form. \newblock {\em Journal
of Dynamic Systems Measurement and Control-Transactions of ASME}, 136(3), 031001, 1-5.

\bibitem{Guo1}
 R. Guo (2014). Stability analysis of a class of switched nonlinear systems with an improved
average dwell time method. \newblock {\em Abstract and Applied Analysis}, 214756, 1-8.


\bibitem{Horn1}
 R. A.  Horn  and  C. R. Johnson (1991). \newblock {\em Topics in Matrix Analysis}, New York: Cambridge University Press.

 \bibitem{Hu1}
 J. Hu and  Y. Hong (2007). Leader-following coordination of multi-agent
systems with coupling time delays.  \newblock {\em Physica A: Statistical Mechanics
and its Applications}, 374(2), 853-863.


 \bibitem{Huang1}   J. Huang (2004). \newblock {\em Nonlinear output regulation: theory and applications}, Phildelphia, PA: SIAM.

\bibitem{Isidori1}
A. Isidori, L. Marconi and G. Casadei (2014). Robust output synchronization of a network of heterogeneous nonlinear agents via nonlinear regulation theory, \newblock {\em IEEE Transactions on Automatic Control}, 59(10), 2680-2691.

\bibitem{Liberzon1}
D.  Liberzon,  (2003), \newblock {\em Switching in Systems and Control}, Boston, MA:
Birkhauser.



\bibitem{Liu2}
 W. Liu  and   J. Huang  (2015). Cooperative global robust output regulation for a class of nonlinear multi-agent systems with switching network,  \newblock {\em IEEE Transactions on Automatic Control}, 60(7), 1963-1968.



%








\bibitem{Praly1}
 L. Praly and  Y. Wang, (1996).  Stabilization in spite of matched unmodelled dynamics and an equivalent definition of input-to-state stability. \newblock {\em Mathematics of Control Signals Systems}, 9(1), 1-33.


%

\bibitem{Sontag1}
  E. D. Sontag and  A. R. Teel, (1995). Changing supply functions in input/state stable systems. \newblock {\em IEEE Transactions on Automatic Control}, 40(8), 1476-1478.





\bibitem{Sontag5}
 E. D. Sontag (1998). Comments on integral variants of ISS.  \newblock {\em Systems \& Control Letters}, 34(1-2), 93-100.







\bibitem{Su4}
Y. Su and J. Huang (2013). Cooperative global output regulation of
heterogenous second-order nonlinear uncertain multi-agent systems.
\newblock {\em Automatica}, 49(11), 3345-3350.


 \bibitem{Su6}
 Y. Su and J. Huang (2015). Cooperative global output regulation for nonlinear uncertain multi-agent systems in lower triangular form.
\newblock {\em IEEE Transactions on Automatic Control}, DOI 10.1109/TAC.2015.2400713.

\bibitem{XuHong1}
D. Xu  and  Y. Hong (2012). Distributed output regulation design for multi-agent systems in output-feedback form. \newblock {\em Proceedings of the 12th International Conference on Control, Automation, Robotics and Vision (ICARCV)}, Guangzhou, China, Dec. 5-7, 596-601.

\bibitem{Zhao1}
J. Zhao, D. J. Hill, and T. Liu (2009).
Synchronization of complex dynamical networks with switching topology: A
switched system point of view. \newblock {\em Automatica}, 45(11), 2502-2511.


%










%
%












%




%




%
%
%






%
%
%
%
%
%
%
%
%
%
%
%
%
%
%
%
%
%
%
%
%
\end{thebibliography}
\end{document}